\documentclass[11pt]{article}
\usepackage[authoryear]{natbib} 

\usepackage{mathrsfs}
\usepackage{graphicx}
\usepackage{mathtools}
\usepackage{caption,subcaption}
\usepackage{amsmath}
\usepackage{amssymb}
\usepackage{textcomp}
\usepackage{euscript}
\usepackage{eufrak}
\usepackage{hyperref}
\usepackage{xcolor}
\usepackage{tcolorbox}
\usepackage[titletoc,title]{appendix}
\usepackage{upgreek}
\usepackage{authblk}

\usepackage{longtable}

\usepackage{float}
\usepackage[scaled]{helvet}

\usepackage{algorithm,algpseudocode}
\algnewcommand{\Inputs}[1]{%
  \State \textbf{Inputs:}
  \Statex \hspace*{\algorithmicindent}\parbox[t]{.8\linewidth}{\raggedright #1}
}
\algnewcommand{\Initialize}[1]{%
  \State \textbf{Initialize:}
  \Statex \hspace*{\algorithmicindent}\parbox[t]{.8\linewidth}{\raggedright #1}
}

\newcommand{\pp}{\mathbb{P}}

\newcommand{\ee}{\mathbb{E}}

\usepackage[utf8]{inputenc}
\usepackage{amsmath} 
\usepackage{amssymb}
\usepackage{amsthm}

\usepackage[margin=1in]{geometry}

\newcommand{\convdistr}{\stackrel{d}{\longrightarrow}}

\newcommand{\lt}{\left}
\newcommand{\rt}{\right}
\newcommand{\wt}{\widetilde}

\newcommand{\var}{\mathbb{V}\mathrm{ar}}

\newcommand{\E}{{\EuScript{E}}} 
\newcommand{\nor}{{\EuScript{N}}} 
\newcommand{\Geo}{{\EuScript{G}}} 

\newcommand{\rmd}{\mathrm{d}}

\newcommand{\whbeta}{\widehat \beta}

\newtheorem{theorem}{Theorem}

\newtheorem{remark}{Remark}

\newtheorem*{algorithm-non}{Algorithm}
\title{Unbiased Estimation of the Reciprocal Mean for Non-negative Random Variables}
\author[*]{Sarat Moka}
\author[*]{Dirk P. Kroese}
\author[**]{Sandeep Juneja}
\affil[*]{\small School of Mathematics and Physics, The University of Queensland, Brisbane}
\affil[**]{School of Technology and Computer Science, Tata Institute of Fundamental Research, Mumbai}
\date{}

\begin{document}

\maketitle
\begin{abstract}
Many simulation problems require the estimation of a ratio of two
expectations. In recent years Monte Carlo estimators have been proposed that can
estimate such ratios without bias.  We investigate the theoretical
properties of such estimators for the estimation of
$\beta = 1/\ee\, Z$, where $Z \geq 0$.  The estimator, $\whbeta(w)$, is of
the form $w/f_w(N) \prod_{i=1}^N (1 - w\, Z_i)$, where $w < 2\beta$ and
$N$ is any random variable with probability mass function $f_w$ on the
positive integers. For a fixed $w$, the optimal choice for $f_w$ is
well understood, but less so the choice of~$w$.  We study the properties of $\whbeta(w)$ as a function
of~$w$ and show that its expected time variance product decreases as
$w$ decreases, even though the cost of constructing the estimator
increases with~$w$. We also show that the estimator is asymptotically
equivalent to the maximum likelihood (biased) ratio estimator and
establish practical confidence intervals.
\end{abstract}

\section{Introduction}
\label{sec:intro}
Over the past few years, unbiased Monte Carlo estimation methods have
received significant attention, due to both theoretical interest and
practical applications; see, for example, \cite{RG15,BCG15,BG15,Jacob15,McLeish2011}. Efficient unbiased estimation of  non-linear
functions of expectations of random variables is  challenging and has
several applications; see, for example, \cite{BCG15,Jacob15}. An
important ``canonical'' case is  the unbiased estimation of $1/\ee\, Z$
for a non-negative random variable~$Z$. Applications include
{regenerative simulation}, estimating parameters involving densities
with unknown normalizing constants, and Bayesian inference. \\

Motivated by these applications, we study properties of an unbiased
estimator of $\beta = 1/\ee\, Z$ proposed by \cite{BCG15} (which is in
turn based on the ideas proposed by \cite{RG15} in
the context of stochastic differential equations). The estimator is
obtained as follows. Write 
 $\beta = \frac{1}{\ee\, Z}  = w\, \sum_{n = 0}^\infty (1 - w\, \ee\, Z)^n$ for $w < 2\beta$; here the condition $w < 2\beta$ guarantees the convergence of the geometric series $ \sum_{n = 0}^\infty (1 - w\, \ee\, Z)^n$.
Further, let $\{Z_i , i \geq 0 \}$ be a sequence of iid copies of $Z$,
and let $N$ be a non-negative integer-valued random variable with $q_n := \pp(N = n) > 0$, for all $n \geq 0$. Then
\begin{align*}
\frac{1}{\ee\, Z} &= w\sum_{n = 0}^\infty q_n \frac{(1 - w\, \ee\, Z)^n}{q_n}= w \sum_{n = 0}^\infty q_n \frac{\ee \prod_{i =1}^n(1 - w\,Z_i)}{q_n} = w\, \ee\lt[\frac{1}{q_{N}} \prod_{i =1}^{N}(1 - w\,Z_i)\rt].
\end{align*}
Define, 
\begin{align}
\label{eqn:ub_estr}
\whbeta(w) := \frac{w}{q_{N}} \prod_{i =1}^{N}(1 - w\,Z_i).
\end{align}
Clearly, ${\ee \whbeta(w) = \beta}$ and thus $\whbeta(w)$ is an unbiased estimator of $\beta$. Note that if $Z \leq b$ almost surely for a constant $b$, then with the choice $w < 1/b$, $\whbeta(w)$ becomes non-negative.  
In this paper, the goal is to study optimal choices for $w$ and ${\{q_n, n  \geq 0\}}$ that make $\whbeta(w)$ {\em efficient}. In particular, a brief description of our contributions is as follows:
\begin{itemize}
\item When $\{q_n, n\geq 0\}$ is the variance-minimizing distribution
  for a fixed $w$, we show that as $w \searrow 0$, the expected cost to construct $\whbeta(w)$ increases to $\infty$, while both the variance and the expected time variance product of $\whbeta(w)$ decease. 
\item As a consequence, we argue that for any $w$, instead of
  approximating $\beta$ with a sample mean of iid copies of
  $\whbeta(w)$, it is optimal to approximate it by just one outcome of
  $\whbeta(w^*)$, where $w^*$ is such that $w^* < w$ and the expected cost of constructing $\whbeta(w^*)$ is the same as that of the sample mean.
\item We study the asymptotic distribution of $\whbeta(w)$ as $w
  \searrow 0$ (i.e., as the expected computational cost for
 the estimator goes to $\infty$). We establish a central limit theorem type convergence result that is useful for finding asymptotically valid confidence intervals.
\item We compare the asymptotic performance of the unbiased estimator $\whbeta(w)$ with that of the maximum likelihood (biased) ratio estimator, where $\beta$ is estimated using the reciprocal of a sample mean of iid copies of $Z$.
\item The above results are studied under the assumption that $N$ has the variance-minimizing distribution. Generating samples from this distribution is impossible as it involves unknown parameters. Since $\whbeta(w)$ is unbiased even for a different distribution of $N$, we develop a method to implement the estimator by proposing a distribution for $N$ (using samples of $Z$) that closely resembles to the variance-minimizing distribution .
\end{itemize}
 
{\em Background:} Several applications of Monte Carlo
simulation involve the estimation of $\beta = 1/\ee\, Z$ for a
non-negative random variable $Z$. In some applications it is a
desirable property to have an unbiased estimator of $\beta$ when the
magnitudes of the available biased estimators are unknown a
priori. Examples include the estimation of a steady-state parameter
$\alpha = \ee R/\ee \tau$ for a {\em regenerative} stochastic process,
where  $\tau $ denotes the length of a regenerative cycle and $R$
denotes the {\em cumulative reward} obtained over the regenerative
cycle; see, e.g, \cite{Glynn06,AG07,MJ14}. It is evident that we have an unbiased estimator of $\alpha$ if we have an unbiased estimator of $1/\ee \tau$. A similar case is where parameters can be expressed as 
$\alpha = {\ee\lt[h(X)f(X)\rt]/\ee f(X)}$ for some real-valued
  function $h$ and probability density $f$, where $f$ is known up to
  the  normalizing constant $\ee f(X)$. Such densities occur, for example,
 in Gibbs point processes (\cite{Moller03}); and a standard method to
 estimate such parameters is by using Markov Chain Monte Carlo (MCMC)
 methods,  see \cite{AG07,RD17}. However, in many situations it is
 difficult to bound the bias of the MCMC estimator, as the mixing time
 of the Markov chain is unknown. An alternative approach is to use a
 ratio estimator, where $\alpha$ is approximated by ratio of the
 sample means of the numerator and the denominator. However, this
 still returns a biased estimator and the bias decreases at a rate
 that is 
 inversely proportional to the sample size; see Remark~\ref{rem:re} and also \cite{AG07}. Therefore, it is desirable to have an unbiased estimator for $1/\ee f(X)$ (and equally for $1/\ee \tau$) that has the same order of complexity as that of the ratio estimator. \\

Most importantly, in some applications, having an unbiased estimator
of $\beta$ is essential. For example, in the study of doubly intractable
models in Bayesian inference, it is assumed that the observations
follow a distribution with a 
density of the form 
$ f(y|\theta) = \frac{g(y,\theta)}{\int g(y, \theta) \rmd y},$
where $g(y, \theta)$ can be evaluated point-wise up to the normalizing
constant $\int g(y, \theta) \rmd y$; see, for example,
\cite{AMYHD15, Walker11,Jacob15}. Standard Metropolis--Hastings algorithms to obtain posterior estimates are not applicable due to the intractability of the normalizing constant. However, an exact inference method called  {\em pseudo-marginal Metropolis--Hastings} proposed by \cite{AR09} can be implemented if a non-negative unbiased estimator of $1/\int g(y, \theta) \rmd y$ is available; also see \cite{Beaumont03,Jacob15,Walker11}. In particular, \cite{Jacob15} highlight the importance of the estimators of the form \eqref{eqn:ub_estr}.\\

A standard method called {\em Russian roulette truncation} can be used for unbiased estimation of~$\beta$. This method is first proposed in the physics literature \cite{CC75,LK91} and further studied by \cite{McLeish2011,GR14, AMYHD15, WM16}. The key drawback of these estimators is that they can take negative values with positive probability even when $Z$ is bounded.\\

{\em Organization of the paper:} In Section~\ref{sec:proper}, we study the properties of the estimator as a function of $w$, under the assumption that $N$ has the variance minimizing estimator. An implementable method is proposed in Section~\ref{sec:impl}.   A conclusion of the paper is given in Section~\ref{sec:conc}. All the results are proved in Appendix~\ref{sec:proofs}.

\section{Properties of the Estimator}
\label{sec:proper}
Without loss of generality, assume that $Z$ is non-degenerate. As the estimator $\whbeta(w)$ in \eqref{eqn:ub_estr} is unbiased, a sample mean of independent copies of
$\whbeta(w)$ is an unbiased estimator of~$\beta$ as well. It is well
known that the sample mean has square-root convergence rate if
$\var\, \whbeta(w) < \infty$; see, e.g., \cite{AG07}. Thus, a
simple strategy is to seek values of $w$ and $\{q_n , n  \geq 0\}$
that minimize $\var\, \whbeta(w)$. Using the  Cauchy--Schwarz
inequality, \cite{BCG15} show that for any $w < 2\,\ee\, Z/\ee\, Z^2$,
$\var\, \whbeta(w)$ is finite and is minimal if $N$ has a  geometric
distribution on the non-negative integers with success probability
$$p_w = 1 - \sqrt{\ee (1 - w\, Z)^2} = 1 - \sqrt{1 - 2\, w\, \ee\, Z
  + w^2\, \ee\, Z^2};$$ that is, if 
\begin{align}
\label{eqn:var_min_pmf}
q_n = (1- p_w)^n p_w, \quad n \geq 0,
\end{align}
where the assumption $w < 2\,\ee\, Z/\ee\, Z^2$ guarantees that $p_w >
0$. Unfortunately, $p_w$ depends on $\ee\, Z$ and $\ee\, Z^2$, which are
unknown. However, $\whbeta(w)$ is unbiased even when $N$ has a
different distribution. Therefore, in the implementation of
$\whbeta(w)$, we can replace $p_w$ with an estimate of $p_w$; see Section~\ref{sec:impl}. In this
section, we study the properties of $\whbeta(w)$ under the assumption
that $N$ has the distribution~\eqref{eqn:var_min_pmf},
because it offers an understanding of what is the best that can be expected from the estimator.\\

Note that the variance of $\whbeta(w)$ is given by,
\begin{align}
\label{eqn:var_Vw}
\var\, \whbeta(w) &= w^2 \sum_{n = 0}^\infty \frac{\lt(\ee(1 - w\, Z)^2\rt)^n}{q_n} - \beta^2 = \frac{w^2}{p_w^2} - \beta^2,
\end{align}
for all $0 < w < 2\,\ee\, Z/\ee\, Z^2$. 
Further, observe that $\ee N = 1/p_w - 1$. Now we can ask what is the value of $w \in (0, 2\,\ee\, Z/\ee\, Z^2)$ that minimizes \eqref{eqn:var_Vw}. This question is not addressed in the existing literature. 
In addition to the variance, it is often important to include the running time to construct the estimator to determine its efficiency; see \cite{GW92}. In that case, we need to select $w$ for which the {\em expected time variance product}, 
$\ee T \, \var\, \whbeta(w),$
is minimal, where $T$ is the time required to construct
$\whbeta(w)$. From \eqref{eqn:ub_estr} (since $Z_1, Z_2, \dots$ are
iid), it is reasonable to assume that $T$ is proportional to the
number of $Z_i$'s used for constructing $\whbeta(w)$. Since $N = N(w)$
samples of $Z$ are used in the construction of $\whbeta(w)$, we assume
that the expected time variance product is $\ee N\, \var\, \whbeta(w)$.

\begin{theorem}
\label{thm:dec_var_g}
Suppose that $N$ has the geometric distribution given in  \eqref{eqn:var_min_pmf}. Then the following hold true.
\begin{itemize}
\item[(i)] The success probability $p_w$ is a strictly concave
  function of $w\in (0, 2\,\ee\, Z/\ee\, Z^2)$ with a  maximum value of $1
  - \sqrt{1 - (\ee\, Z)^2/ \ee\, Z^2}$ attained at $w = \ee\, Z/\ee\, Z^2$, and 
$$ \lim_{w \searrow 0} p_w = \lim_{w \nearrow 2\,\ee\, Z/\ee\, Z^2} p_w = 0.$$
\item[(ii)] The variance $\var\, \whbeta(w)$ is a strictly increasing
  convex function of $w \in (0, 2\,\ee\, Z/\ee\, Z^2)$,  with
	\begin{align*}
 	\lim_{w \searrow 0}\var\, \whbeta(w)~=~ 0 \quad \text{ and } \quad  	\lim_{w \nearrow 2\,\ee\, Z/\ee\, Z^2}\var\, \whbeta(w)~=~ \infty. 
	\end{align*} 
\item[(iii)] The expected time variance product ${\ee N(w)\, \var\, \whbeta(w)}$ is a
  strictly increasing function of $w \in (0, 2\,\ee\, Z/\ee\, Z^2)$, with
\begin{align*} 
\lim_{w \searrow 0} \ee N(w)\, \var\, \whbeta(w) = \frac{\var\, Z}{(\ee\, Z)^4}, \quad \text{ and } \quad \lim_{w \nearrow 2\,\ee\, Z/\ee\, Z^2} \ee N(w)\, \var\, \whbeta(w) = \infty.
\end{align*}
\end{itemize}
\end{theorem}

\begin{remark}
\label{rem:opt_w}
\normalfont
To understand the implications of Theorem~\ref{thm:dec_var_g}, suppose
we select a $w \in \lt(0, 2\,\ee\, Z/\ee\, Z^2\rt)$, giving an expected computational
cost of $\ee
N(w)$ to obtain $\whbeta(w)$. Further, let $\overline
\beta_k(w)$ be the sample mean of $k$ iid copies of
$\whbeta(w)$. Then the expected time variance product for $\overline \beta_k(w) $ is 
 \[
 k\, \ee N(w)\, \var\, \overline \beta_k(w)  = \ee N(w)\, \var\, \whbeta(w).
 \]
 Now suppose that $w^* < w$ is selected so that the average cost to
 generate one outcome of $\whbeta(w^*)$ is equal to the average cost to
 construct $\overline \beta_k(w)$; that is, $k\, \ee N(w)$. Then, from
 Theorem~\ref{thm:dec_var_g} (iii), for the same computational effort,
 $\widehat{\beta}(w^*)$ has a
 smaller variance than $\overline\beta_k(w)$ for any feasible $w$
 selected as above,  and is therefore a better estimator.
\begin{figure*}[h!]
\hspace{-5mm}
    \centering
    \begin{subfigure}[t]{0.45\textwidth}
        \centering
        \includegraphics[height=0.7\textwidth, width=1.1\textwidth]{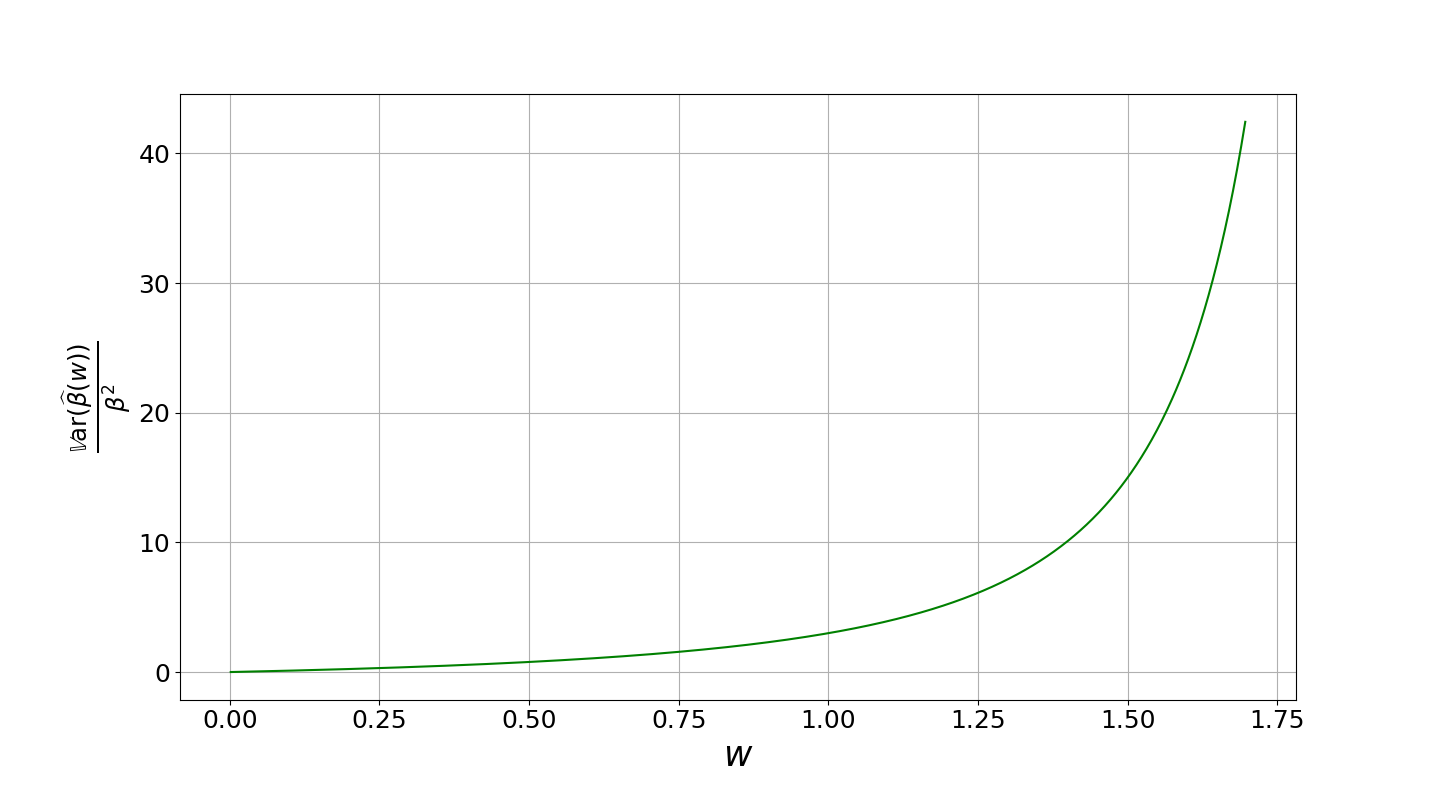}
        \caption{ }
    \end{subfigure}%
    \hspace{7mm}
    \begin{subfigure}[t]{0.45\textwidth}
        \centering
        \includegraphics[height=0.7\textwidth, width=1.1\textwidth]{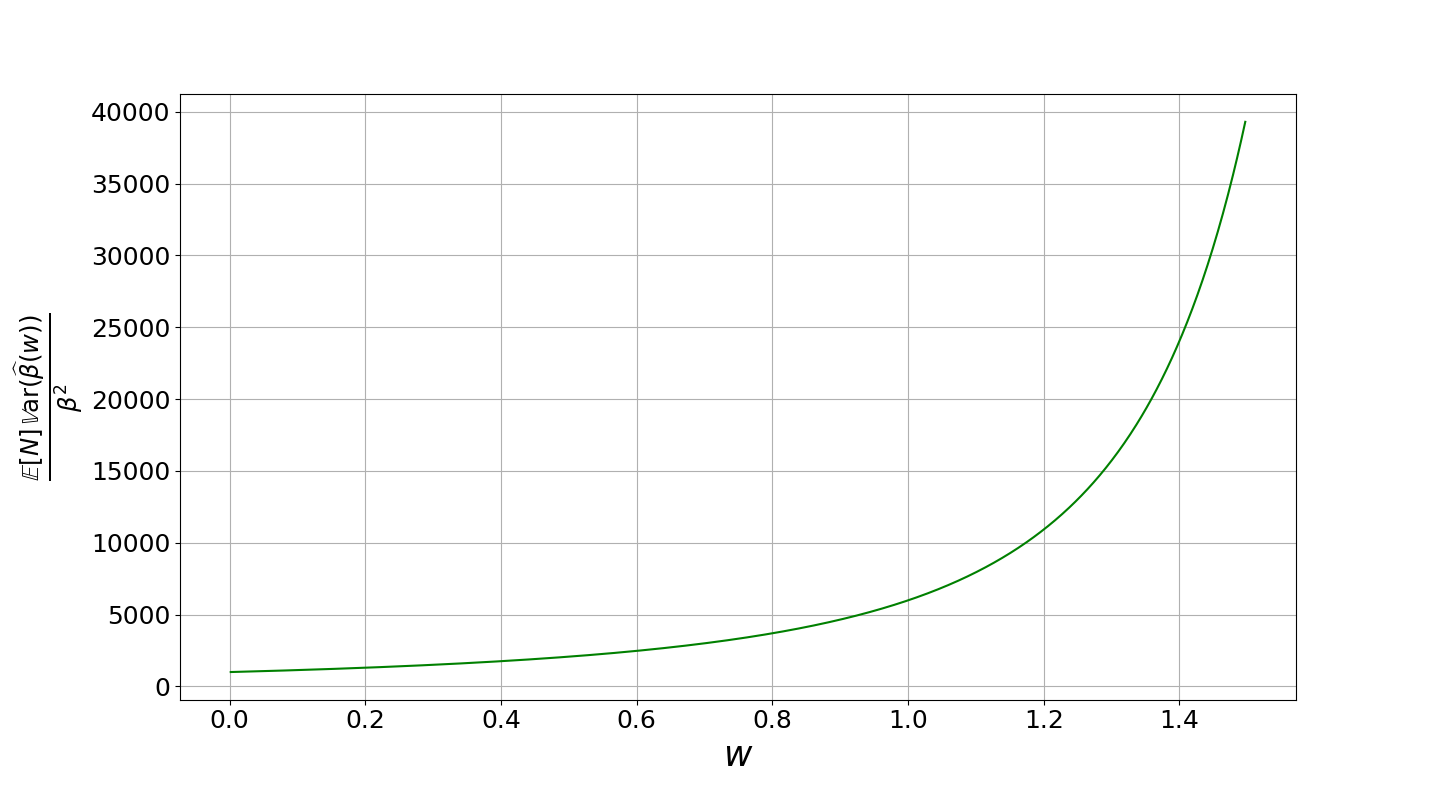}
        \caption{}        
    \end{subfigure}%
    \caption{\footnotesize An example to illustrate the dependency  of the performance of the unbiased estimator on parameter~$w$. Panels (a)  and (b) show, respectively, the relative variance $\frac{\var\, \whbeta(w)}{\beta^2}$ and the expected time relative variance product $\frac{\ee N(w)\var\, \whbeta(w)}{\beta^2}$, as functions of~$w$. }
    \label{fig:w_choice}
\end{figure*}

To illustrate the results of Theorem~\ref{thm:dec_var_g}, consider an example where
${Z = \mathbb{I}(A)}$ for an event $A$ with probability ${\pp(A) = \ee\, Z = 0.001}$. Since ${\ee\, Z^2 = 0.001}$, the relative variance ${\frac{\var\, Z}{(\ee\, Z)^2} = 999}$. 
By substituting the values of $\ee\, Z$ and $\ee\, Z^2$, we can calculate $p_w$, $\ee N $ and $\var\, \whbeta(w)$ for any ${w < 2\,\ee\, Z/\ee\, Z^2 = 2}$.  Figure~\ref{fig:w_choice} illustrate the effect of $w$ on the efficacy of the estimator~$\whbeta(w)$. As expected, both $\frac{\var\, \whbeta(w)}{\beta^2}$ and $\frac{\ee N(w)\var\, \whbeta(w)}{\beta^2}$ are decreasing as $w \searrow 0$ with the limits~$0$ and $ {\frac{\var\, Z}{(\ee\, Z)^2} = 999}$, respectively.
\qed
\end{remark} 

Remark~\ref{rem:opt_w} motivates us to study the asymptotic
distributional properties of $\whbeta(w)$ as $w \searrow 0$, when $N$
has the geometric distribution given in \eqref{eqn:var_min_pmf}.
Theorem~\ref{thm:conv_V} is crucial for establishing confidence intervals that are asymptotically valid as $w \searrow 0$. 

 \begin{theorem}
 \label{thm:conv_V}
Suppose that $N$ has the distribution given by \eqref{eqn:var_min_pmf} and $Z$ is bounded. Then, as $w \searrow 0$,
\begin{align*}
(i)\quad \whbeta(w)  \convdistr \beta,\quad \quad\text{ and }\quad \quad (ii)\quad \frac{\lt( \whbeta(w) - \beta \rt)}{\sqrt{w\, \ee\, Z}}  \convdistr \lt(\sqrt{\frac{\var\, Z}{(\ee\, Z)^4} \E(1)}\rt)\nor(0,1),
\end{align*}
where  $\convdistr$ denotes  convergence in distribution, and
$\E(1)$ and $\nor(0,1)$ are
independent  random variables from respectively the
standard (mean 1) exponential and standard normal distributions.\\
 \end{theorem}
 
We show later in Section~\ref{sec:proof_thm1}  that $\frac{p_w}{w\, \ee\, Z} \to 1$ as $w \searrow 0$ (see \eqref{eqn:pw_bnd}). Since $\ee N = 1/p_w - 1$ and $\lim_{w \to 0}p_w = 0$,
$w\, \ee Z\: \ee N \to 1$ as $w \searrow 0$. That means, an alternative expression for Theorem \ref{thm:conv_V} (ii) is 
\[
\sqrt{\ee N}\lt( \whbeta(w) - \beta \rt)  \convdistr \sigma\sqrt{
  \E(1)}\, \nor(0,1),\quad \text{ as} \quad w \searrow 0,
\]
where $\sigma = \sqrt{\frac{\var\, Z}{(\ee\, Z)^4}}$. The above expression has more resemblance to the standard central limit theorem, since $\ee N$ is the computational cost of the estimator.
It is not difficult to show that $\sqrt{\E(1)} \, \nor(0,1)$ is a random variable with density 
$f(x) = \frac{1}{\sqrt{2}} \exp(- \sqrt{2} |x|)$, which is the density of a mean zero {\em Laplace} (or double exponential) distribution with scale $1/\sqrt{2}$.
These observations are useful for constructing asymptotically valid
confidence intervals as follows. For any given $\alpha \in (0,1)$, by
solving $\int_{0}^t f(x) \,dx = (1 - \alpha)/2$ for $t$, we get $t = -\log(\alpha)/\sqrt{2}$.
Then using Theorem~\ref{thm:conv_V}, we can say that the interval $$\lt(\whbeta(w) + \frac{\log(\alpha)}{\sqrt{2}} \sigma \sqrt{w\, \ee\, Z}, \quad \whbeta(w) -\frac{\log(\alpha)}{\sqrt{2}} \sigma \sqrt{w\, \ee\, Z}\rt)$$ is an asymptotic $1 - \alpha$ confidence interval for $\beta$.

\begin{remark}[Comparison with the ratio estimator]
\label{rem:re}
\normalfont
A standard (biased) estimator of $\beta$ is ${1/ \bar Z_n}$, where
$\bar Z_n$ is the sample mean of $n$ iid copies of $Z$. Using Taylor's
theorem for the function $1/x$ about $\ee\, Z$, we can easily show that
the bias of $1/\bar Z_n$ is approximately $\frac{1}{n}\frac{\var\, Z}{(\ee\, Z)^3}$ for large $n$, while, on the other hand,
$\whbeta(w)$ has zero bias. Furthermore, using the same Taylor's theorem, we can show that  the
asymptotic time variance product of $1/\bar Z_n$ is
$\var\, Z/(\ee\, Z)^4$ as $n \to \infty$. From
Theorem~\ref{thm:dec_var_g} (iii), the unbiased estimator $\whbeta(w)$ has the same asymptotic expected time variance
product. However, unbiasedness of $\whbeta(w)$ comes at cost. As an application of the delta method, we can show that the ratio estimator satisfies the central limit theorem:
$\sqrt{n} \lt(1/ \bar Z_n - \beta \rt) \convdistr \sqrt{\frac{\var\,
    Z}{(\ee\, Z)^4}} \, \nor(0,1).$ That is, the ratio estimator is asymptotically normal.
On the other hand, the asymptotic distribution of the unbiased estimator $\whbeta(w)$ is Laplace, which has more slowly decaying tails than a normal distribution. In conclusion, the ratio estimator can have narrower confidence intervals than the unbiased estimator.
\end{remark}

 \begin{remark}[Importance sampling]
 \label{rem:imp_sam}
 \normalfont
Just like in the case of the ratio estimator, from Theorems
\ref{thm:dec_var_g} and \ref{thm:conv_V}, the relative variance of $Z$
is the key factor influencing the asymptotic properties of the
unbiased estimator. The smaller the value is of the relative variance
of $Z$, the better is the reliability of the unbiased estimator. One of the most effective technique of variance reduction is {\em importance sampling}. We can improve the performance of the estimator if we can implement an importance sampling technique  for  the random variable~$Z$.
 \end{remark}

\begin{remark}[The time variance product minimizing distribution for $N(w)$]
\label{rem:tvmd}
\normalfont
We have assumed that for a given $w$ the random variable $N(w)$ has the
variance minimizing distribution given by
\eqref{eqn:var_min_pmf}. However, when the criteria for the optimality
of $\whbeta(w)$ is the  minimization of the expected time variance
product, we need to seek a distribution $\{q_n, n \geq 0\}$
that minimizes $\ee N(w)\var\, \whbeta(w)$. It is shown in
\cite{BCG15} 
that the
distribution that minimizes the expected
time variance is given by 
\begin{align}
\label{eqn:time_var_min_pmf}
\wt q_n = \frac{w\, (1- p_w)^n}{\sqrt{\beta^2 + d_w\, n}}, n \geq 0,
\end{align}
where $d_w$ is the unique (positive) number satisfying $\sum_{n = 0}^\infty \frac{w\,(1- p_w)^n}{\sqrt{\beta^2 + d_w\, n}} = 1.$
When compared to the distribution \eqref{eqn:var_min_pmf}, drawing
samples from \eqref{eqn:time_var_min_pmf} has an extra difficulty of
finding $d_w$ by solving an equation that contains an infinite
sum. Even if we overcome this difficulty, the reduction in the
expected time variance product is typically insignificant, because for small values of $w$,
\begin{align}
\label{eqn:tvp_order}
\ee \wt N(w)\, \var\, \wt \beta(w) = \ee N(w)\, \var\, \whbeta(w) (1  + O(w)),
\end{align}
when  $N(w)$ has the variance minimizing distribution~\eqref{eqn:var_min_pmf}; see Section~\ref{sec:tvp_order_proof} for a proof of \eqref{eqn:tvp_order}.
\qed
\end{remark}

\section{An Implementation}
\label{sec:impl}
Recall that the success probability $p_w$ of $N$ is a function of unknown quantities $\ee\, Z$ and $\ee\, Z^2$. However, fortunately, $\whbeta(w)$ in \eqref{eqn:ub_estr} is still an unbiased estimator of $\beta$ for any distribution $\{q_n, n \geq 0\}$ of $N$. In particular, instead of taking $q_n$ as in \eqref{eqn:var_min_pmf}, we can take $q_n = P_k(1 - P_k)^n$, where $P_k$ is defined below. Under the proposed implementation, when the given budget is sufficiently large, half of the budget is used for estimating $P_k$ and then $w$ is chosen such that the remaining half the budget is used for generating a sample of the unbiased estimator.\\

To simplify the discussion, assume that there is a known constant $0 <
\varepsilon \leq \beta$; for example, if $Z \leq b$ for a constant $b$, we can take $\varepsilon = 1/b$. Let $\wt Z_1, \wt Z_2, \dots, \wt Z_k$ be a sequence of iid copies of $Z$, independent of the sequence $Z_1, Z_2, \dots$, which is used in the construction of the unbiased estimator $\whbeta$ in \eqref{eqn:ub_estr}. 
Define the first two sample moments: ${M_1(k) = \frac{1}{k}\sum_{i = 1}^k \wt Z_i}$ and ${M_2(k) = \frac{1}{k}\sum_{i = 1}^k \wt Z^2_i}$. If $M_1(k) > 0$, define, 
$$P_k = 1 - \sqrt{ \frac{1}{k} \sum_{i = 1}^k (1 - w_k\, \wt Z_i)^2} \quad \text{with} \quad w_k = \min\lt(\frac{1}{k\,M_1(k)}, \frac{M_1(k)}{M_2(k)}, \varepsilon\rt).$$
Otherwise, take $P_k = 1/k$ and $w_k = \varepsilon/k$. The condition $
w_k < 2 \frac{M_1(k)}{M_2(k)}$ guarantees that $P_k > 0$. Further,
whether $M_1(k) = 0$ or not, $w_k  < \beta$ and hence it guarantees
that the estmator $\whbeta(w_k)$ (defined by \eqref{eqn:ub_estr}) with
$q_n =P_k(1 - P_k)^n $ is an unbiased estimator of~$\beta$. Note that
given $M_1(k)$ and $M_2(k)$,
the expected cost to construct to $\whbeta(w_k)$ is $k + 1/P_k$
(including the cost to construct $P_k$),  since $\ee\lt[N(w_k) \big|M_1(k), M_2(k)\rt] = 1/P_k - 1$. Theorem~\ref{thm:impl} states that for large values of $k$, this total expected cost is approximately $2\,k$, and conditioned on $M_1(k)$ and $M_2(k)$, the expected time variance product goes to a random variable with mean $4\frac{\var\, Z}{(\ee\, Z)^4}$. See Section~\ref{sec:proof_impl} for a proof Theorem~\ref{thm:impl}. 
\begin{theorem}
\label{thm:impl}
 Under the above construction, $\ee \whbeta(w_k) = \beta$, and as $k \rightarrow \infty$, $k\, P_k \to 1,\, a.s.$, and 
$$\lt(k + \frac{1}{P_k}\rt)\var\,\lt(\whbeta(w_k) \big| M_1(k), M_2(k)\rt)\longrightarrow 2\frac{\var\, Z}{(\ee\, Z)^4}  \lt[1 + \chi^2_1 \rt],\, a.s.,$$
where $ \chi_1^2$ is the square of a standard normal random variable.
\end{theorem}
To understand  Theorem~\ref{thm:impl}, consider the example given in Remark~\ref{rem:opt_w}. We estimated the expected total cost $k + \ee[1/P_k]$ and $\var\, \whbeta(w_k)$ using $10000$ samples of $P_k$ and $\whbeta(w_k)$, respectivley, with $k = 10000$. Our simulation results show that the estimated expected time relative variance product is $3969.75$, which is approximately equal to $4 \frac{\var\, Z}{(\ee\, Z)^2} = 4\times 999 = 3996$, as expected.

\section{Conclusion}
\label{sec:conc}
We investigated the theoretical properties of a parametrized family
$\{ \whbeta(w) , w > 0\}$ of unbiased estimators of $1/\ee\, Z$
for a non-negative random variable $Z$. We studied the variance and
the expected time variance product as functions of $w$ and established
several asymptotic results. We showed that with an optimal choice of
$w$, the asymptotic performance of the unbiased estimator $\whbeta(w)$
is comparable to the maximum likelihood (biased) ratio estimator. We
further proposed an implementable unbiased estimation based on our
results. Similar to Theorem~\ref{thm:conv_V}, our ongoing research
establishes a central limit theorem type convergence result for
$\whbeta(w_k)$ defined in Section~\ref{sec:impl}, by taking the budget parameter $k \to \infty$.


\appendix
\section{Appendix}
\label{sec:proofs}
To simplify the notation in this section, we use $z_1 :=
\ee\, Z$ and $z_2 := \ee\, Z^2$. We also write $p_w'$ for the derivative
$\frac{\rmd p_w}{\rmd w}$ and $p_w''$ for the second derivative. $\E(\lambda)$ and $\Geo(p)$ denote respectively the mean $1/\lambda$ exponential distribution and the geometric distribution on non-negative integers and the success probability~$p$. 

\subsection{Proof of Theorem~\ref{thm:dec_var_g}}
\label{sec:proof_thm1}
From the definition, $p'_w  =  (z_1 -w\, z_2)/(1-p_w).$
It follows that $p''_w = - \frac{1}{(1 - p_w)^3}\lt[z_2 - z_1^2\rt] < 0,$
where the strict inequality holds because $z_2 > z_1^2$, which follows from the assumption that $Z$ is non-degenerative. Therefore, $p_w$ is strictly concave over $(0, 2z_1/z_2)$ and it achieves its maximum value $1 - \sqrt{1 - z_1^2/z_2}$ at $w = z_1/z_2$. From the definition of $p_w$, it is evident that
$\lim_{w \searrow 0} p_w = \lim_{w \nearrow 2z_1/z_2} p_w = 0$.

Recall from \eqref{eqn:var_Vw} that the variance of the estimator is
equal to ${w^2/p_w^2} - z_1^2$. Its derivative can be written as
\begin{align}
\label{eqn:der_var_V}
\frac{\rmd \var\, \whbeta(w)}{\rmd w} = 2w \, \frac{w\, z_1 - p_w }{p_w^3(1 - p_w)} \end{align}
and the second derivative as
\begin{equation}
\frac{\rmd^2 \var\, \whbeta(w)}{\rmd w^2} = \frac{2(w\, z_1 - p_w)}{p_w^3(1 - p_w)^3}\lt(3\lt(w\, z_1 - p_w\rt) + 2p_w \lt( p_w - w^2\, z_2\rt)\rt).
\end{equation}

Using Jensen's inequality, ${\ee (1 - w\, Z)^2 >  (1 - w\, z_1)^2},$
where the strict inequality holds again because $Z$ is
non-degenerative. On the other hand, by  Bernoulli's inequality, $$\sqrt{\ee (1 - w\, Z)^2 }  = \sqrt{ 1 + (-2w\, z_1 + w^2\, z_2)}$$ is maximized by
$1 - w\, z_1 + w^2\, z_2/2.$
Thus,
\begin{align}
\label{eqn:pw_bnd}
z_1 - w\, z_2/2 \leq \dfrac{p_w}{w} < z_1.
\end{align}
Using \eqref{eqn:pw_bnd}, we have $w\, z_1 - p_w > 0$ and $p_w - w^2\, z_2 \geq 2wz_1 - p_w > w\, z_1$, and hence for all $w \in (0, 2z_1/z_2)$, $\frac{\rmd \var\, \whbeta(w)}{\rmd w} > 0$ and $\frac{\rmd^2 \var\, \whbeta(w)}{\rmd w^2} > 0,$
which  establishes the convexity of $\var\, \whbeta(w)$ over $(0, 2z_1/z_2)$.\\

We now prove that the expected time variance product is a strictly increasing over $(0, 2z_1/z_2)$. 
Let $g_1(w) := \ee N(w)$, $g_2(w) := \frac{\var\, \whbeta(w) }{\beta^2}$ and $g(w) := g_1(w)g_2(w) = \lt(\frac{1}{p_w} -1 \rt)\lt(\frac{w^2\, z_1^2}{p_w^2} - 1\rt).$
Then,
\begin{align*}
\frac{\rmd g_1}{\rmd w}(w) = \frac{w\, z_2 - z_1}{p_w^2(1 - p_w)},\quad \text{ and } \quad \frac{\rmd g_2}{\rmd w}(w) = 2w\, z_1^2 \lt[ \frac{w\, z_1- p_w}{p_w^3(1 - p_w)}\rt],
\end{align*}
and hence, we write
\begin{align*}
\frac{\rmd g}{\rmd w}(w) &= \frac{2w\, z_1^2(1 - p_w)(w\, z_1 - p_w) + (w^2\, z_1^2 - p_w^2)(w\, z_2 - z_1)}{p_w^4(1 - p_w)}\\
                         &= \frac{wz_1 - p_w}{p_w^4(1 - p_w)} \lt[ w\, p_w(z_2 -  z_1^2) + z_1 (w\, z_1 - p_w) + w\, z_1(w\, z_2 - z_1p_w)\rt] >0,
\end{align*}
where the inequality holds because $z_2 > z_1^2$, $w\, z_1 > p_w$ and $w\, z_2 > w\, z_1^2 > z_1\, p_w$. Therefore, $\ee N(w)\, \var\, \whbeta(w)$ is strictly increasing over $(0, 2z_1/z_2)$.\\

The claims that $\lim_{w \nearrow 2z_1/z_2} \var\, \whbeta(w) = \infty$ and $\lim_{w \nearrow 2z_1/z_2} \ee N(w)\var\, \whbeta(w) = \infty$ hold trivially because $\lim_{w \nearrow 2z_1/z_2} p_w = 0$.
To complete the proof of the theorem, we can write, by Taylor's
theorem, for any $x \in (0,1)$: 
$\sqrt{x} = 1 + \frac{(x - 1)}{2} - \frac{(x - 1)^2}{8} + \frac{(x -1)^3}{16 \tilde{x}^{5/2}},$
for some $\tilde x \in (x, 1).$ Consequently,
\begin{align*}
\sqrt{\ee (1 - w\, Z)^2} = 1 + \frac{\ee (1 - w\, Z)^2 - 1 }{2} - \frac{\lt(\ee (1 - w\, Z)^2 - 1\rt)^2}{8} + R(w),
\end{align*}
where $R(w) = \frac{\lt(\ee (1 - w\, Z)^2 - 1 \rt)^3}{16 \tilde x^5/3}$ for some $\tilde x \in \lt(\ee (1 - w\, Z)^2, 1\rt)$. Since $\tilde x \to 1$ as $w \to 0$ and $ \ee (1 - w\, Z)^2 - 1 = w^2 z_2 - 2w\, z_1$, we have $
R(w) = O(w^3)$. Further simplification yields that $p_w =  w\, z_1 - \frac{w^2}{2}\lt(z_2 - z_1^2\rt) + O(w^3),$ and thus,
$\frac{p^2_w}{w^2\, z^2_1} = 1 - \frac{w}{z_1}\lt(z_2 - z_1^2\rt) + O(w^2).$
Since $1 - p_w = 1 + O(w)$, 
\begin{align}
\label{eqn:cost_var}
\frac{1 - p_w}{p_w}  =  \frac{1}{wz_1}\lt( 1 + O(w) \rt), \quad \text{ and } \quad \frac{w^2 z_1^2}{p^2_w} - 1 = w \frac{(z_2 - z_1^2)}{z_1}\lt(1 + O(w)\rt).
\end{align}
We conclude that both $\var\, \whbeta(w)$ and $\ee N(w)\var\, \whbeta(w)$ go to their respective minima as $w \searrow 0$.

\subsection{Proof of Theorem~\ref{thm:conv_V}}
Statement (i) follows directly from Theorem~\ref{thm:dec_var_g} and
Chebyshev's inequality:
$$\pp\lt(|\whbeta(w) - \beta| > \epsilon \rt) \leq
\var\,\whbeta(w)/\epsilon^2 \to 0,\, \text{ as }\, w \searrow 0,$$ for
every $\epsilon > 0$. To prove (ii), consider a decreasing sequence ${
  w_1 > w_2 > \cdots} $ such that $w_1 \leq z_1/z_2$ and $\lim_{k \to
  \infty} w_k = 0$. We construct an almost surely increasing sequence
$N_1 \leq N_2 \leq \cdots$ such that $N_k \sim \Geo(p_{w_k})$
and 
\begin{align}
\label{eqn:finite_wN}
\lim_{k \to \infty} [w_kN_k] = X_\infty/z_1, \quad a.s.,
\end{align} 
for a random variable $X_\infty \sim  \E(1)$. To do so we invoke Theorem 3.1 of \cite{MJ14}. Let $\lambda_k = - \log(1 - p_{w_k})$ and $E_k \sim \E(\lambda_k)$. Then, \cite{MJ14} says that for each $k$, there exist a random variable $Y_k$ with cumulative distribution function 
$$G_k(x) = 1 - \lt(1 - \frac{\lambda_{k+1}}{\lambda_k}\rt) \exp(-\lambda_{k+1} x),\quad  x \geq 0,$$
 such that $Y_k$ is independent of $E_k$,  and $E_{k+1}$ has the same  distribution as $E_k + Y_k$. Therefore, without loss of generality we assume that there is sequence of independent random variables $Y_k \sim G_k(x)$
such that $E_{k+1} = E_{k} + Y_k = E_1 + \sum_{i = 1}^k Y_i$  for all $k \geq 1$.\\

Consider the natural filtration $\lt\{\mathscr{F}_k = \sigma(E_1, \dots, E_k), k \geq 0\rt\}$. Since $E_{k + 1} = E_k + Y_k$,
\begin{align*}
\lambda_{k+1}  E_{k+1} - 1 &= \lambda_{k+1} [ E_{k}  + Y_k - 1/\lambda_k - \ee Y_k]
                                                                   = \frac{\lambda_{k+1}}{\lambda_k} \lambda_k [ E_{k}   -1/\lambda_k] + \lambda_{k+1} [Y_k - \ee Y_k]\\
                                                                   &\leq \lambda_k E_{k}   - 1+ \lambda_{k+1} [Y_k - \ee Y_k],
\end{align*}
where the last inequality holds because $\lambda_{k +1} \leq \lambda_k$. We have $\ee\lt[ \lambda_{k+1} E_{k+1} | \mathscr{F}_k \rt] \leq \lambda_kE_{k}$ since $Y_k$ is independent of $\mathscr{F}_k$.
Thus, $\displaystyle \lt\{ X_k := \lambda_{k} E_k, k \geq 1 \rt\}$ is a supermartingale (with respect to $\lt\{\mathscr{F}_k\rt\}$). 
In fact the sequence $\lt\{X_k, k \geq 1 \rt\}$ is bounded in $\mathcal{L}^2$, because $\sup_{k \geq 1} \ee X_k^2  = 2$, making it uniformly integrable submartingale. Thus, ${X_\infty = \lim_{k \to \infty} X_k}$ exists $a.s.$ (see Theorem~2 in Section~4 of Chapter~VII of \cite{SHI96}).
Since $\pp\lt(X_k \leq x \rt) = \pp\lt(E_k \leq \frac{x}{\lambda_k}  \rt) = 1 - \exp\lt( - x \rt).$
This implies that $X_\infty \sim \E(1)$.\\

 Let $N_k = \lfloor E_k \rfloor$. Then for all $k$ we have $N_k  \sim \Geo(p_{w_k})$ and $N_k \leq N_{k+1}$.
From \eqref{eqn:pw_bnd}, $\lim_{k \to \infty} (1 - p_{w_k})^{1/w_k} = \exp(-z_1)$ and hence
$\lim_{k \to \infty} \frac{w_k}{\lambda_k} = 1/z_1$.  
From the convergence of the sequence $X_1, X_2, \dots$, we have $ \frac{w_k}{\lambda_k}X_k - w_k \leq w_k N_k \leq \frac{w_k}{\lambda_k} X_k$. Therefore, \eqref{eqn:finite_wN} holds.

Now define 
\begin{align}
\label{eqn:V_k}
\whbeta_k := \frac{w_k}{(1 - p_{w_k})^{N_k} p_{w_k}} \prod_{i =1}^{N_k}(1 - w_kZ_i).
\end{align}
From the definitions, $\whbeta_k$ is identical to $\whbeta(w_k)$ in distribution. 
We now conclude the proof Theorem~\ref{thm:conv_V} by establishing lower and upper bounds on $\whbeta_k$ separately. Let $b$ be an upper bound on $Z$. From the construction of $\whbeta_k$ given by \eqref{eqn:V_k}, for all $k$ such that $w_k < 1/b$, we have using  \eqref{eqn:pw_bnd} that 
\begin{align*}
\whbeta_k &\geq \frac{1}{z_1(1 - p_{w_k})^{N_k}} \exp\lt( \sum_{i = 1}^{N_k} \log(1- w_kZ_i)\rt)\geq \frac{1}{z_1} \exp\lt( N_k w_k (z_1 - w_k z_2/2) + \sum_{i = 1}^{N_k} \log(1- w_kZ_i)\rt).
\end{align*}
Using Taylor's theorem,  $\log(1 - x)  \geq -x -  \frac{x^2}{2(1 - x)^2}$ for any $x \geq 0$,  and thus,
\begin{align}
\whbeta_k &\geq \frac{1}{z_1} \exp\lt( N_k w_k z_1 - N \frac{w_k^2 z_2}{2}  - w_k\sum_{i = 1}^{N_k} Z_i -  \sum_{i = 1}^{N_k} \frac{w_k^2Z^2_i}{2(1 - w_kZ_i)^2}\rt)\nonumber\\
        &= \frac{1}{z_1} \exp\lt( - N_k w_k \frac{1}{N_k}\sum_{i = 1}^{N_k} (Z_i - z_1) - \frac{N_kw_k^2 }{2}\lt( z_2  + \frac{b^2}{(1- w_kb)^2}\rt)\rt).\label{eqn:lower_V}
\end{align}
On the other hand, from \eqref{eqn:V_k} and \eqref{eqn:pw_bnd},
\begin{align}
\whbeta_k &\leq \frac{w_k}{(1 - p_{w_k})^{N_k}p_{w_k}} \exp\lt( -\sum_{i = 1}^{N_k} w_k Z_i \rt) \leq \frac{w_k}{(1 - w_k z_1)^{N_k}p_{w_k}} \exp\lt( -\sum_{i = 1}^{N_k} w_k Z_i \rt)\nonumber\\
    &= \frac{w_k}{p_{w_k}} \exp\lt( -\sum_{i = 1}^{N_k} w_k Z_i - N_k \log(1 - w_k z_1)\rt)\nonumber\\
     &= \frac{w_k}{p_{w_k}} \exp\lt( -w_kN_k \frac{1}{N_k}\sum_{i = 1}^{N_k}  (Z_i - z_1) + N_k w^2_k \frac{z_1^2}{2(1 - w_kz_1)^2}\rt).\label{eqn:upper_V}
\end{align}
Using the strong law of large numbers and Theorem~1 of \cite{Rich65}, 
$$\lim_{k \to \infty} \frac{1}{N_k}\sum_{i = 1}^{N_k}  (Z_i - z_1) = 0, \, a.s.$$
Further, using \eqref{eqn:pw_bnd} and
\eqref{eqn:finite_wN}, we have that 
$ \lim_{k \to \infty } \whbeta_k = \beta$ almost surely. From Taylor's
theorem with a Cauchy remainder term, we have almost surely
\begin{align*}
\log(\whbeta_k)  - \log(\beta)  &=  \frac{(\whbeta_k - \beta)}{\beta} - \frac{(\whbeta_k - \widehat X) (\whbeta_k - \beta)}{\widehat X^2}= \frac{(\whbeta_k - \beta)}{\beta} \lt[1 + o(1)\rt]
\end{align*}
for a random variable $\widehat X$ that takes values between $\whbeta_k$ and $\beta$. Therefore, to complete the proof of (ii), it is enough to show that
$
\frac{1}{\sqrt{w_k}}\log\whbeta_k \convdistr \lt(\sqrt{\var\, Z\, X_\infty}\rt)\nor(0,1).
$
From \eqref{eqn:lower_V},
\begin{align}
\label{eqn:upper_clt}
\frac{1}{\sqrt{w_k}}\log[z_1\whbeta_k]  \geq - \sqrt{w_k N_k} \frac{1}{\sqrt{N_k}}\sum_{i = 1}^{N_k} (Z_i - z_1) - \frac{N_kw_k^{3/2} }{2}\lt( z_2  + \frac{b^2}{(1- w_kb)^2}\rt).
\end{align}
From \eqref{eqn:finite_wN} and because $w_k \searrow 0$, the second term on the right hand side of the expression goes to zero. Again using \eqref{eqn:finite_wN}, we conclude that the right hand side of \eqref{eqn:upper_clt} goes to $ \lt(\sqrt{\frac{\var\, Z}{z_1}X_\infty}\rt) \nor(0,1)$ in distribution. 
From \eqref{eqn:upper_V},
\begin{align*}
\frac{1}{\sqrt{w_k}}\log[z_1\whbeta_k]  \leq \frac{1}{\sqrt{w_k}} \log\lt[\frac{z_1 w_k}{p_{w_k}}\rt] - \sqrt{w_k N_k} \frac{1}{\sqrt{N_k}}\sum_{i = 1}^{N_k} (Z_i - z_1) + N_kw_k^{3/2} \frac{ z_1^2}{2(1- w_kz_1)^2}.
\end{align*}
We complete the proof because from \eqref{eqn:pw_bnd}, as $k \to \infty$,
\[
0 \leq \frac{1}{\sqrt{w_k}} \log\lt[\frac{z_1 w_k}{p_{w_k}}\rt] \leq \frac{1}{\sqrt{w_k}} \log\lt[\frac{1}{\lt(1 - \frac{w_kz_2}{2z_1}\rt)}\rt] = - \sqrt{w_k} \log\lt[\lt(1 - \frac{w_kz_2}{2z_1}\rt)^{1/w_k}\rt] \to 0, \quad a.s.
\]

\subsection{Proof of \eqref{eqn:tvp_order}}
\label{sec:tvp_order_proof}
First, observe from the definitions that $\ee \wt N(w)\var\, \wt \beta(w) \leq \ee N(w)\var\, \whbeta(w).$
Further using the fact that $\beta = 1/z_1$ and \eqref{eqn:pw_bnd},
\begin{align*}
\ee \wt N(w) &\geq  \sum_{n=1}^\infty n \frac{p_w(1-p_w)^n}{\sqrt{1 + d_w m^2_1 n}}\\
               &= \ee N(w) \sum_{n=1}^\infty \frac{1}{\sqrt{1 + d_w z_1^2 n}} \frac{ n p_w(1-p_w)^{n}}{\ee N(w)}\\
               &\geq  \ee N(w)  \frac{1}{\sqrt{1 + d_w z_1^2 \ee N(w)^2/\ee N(w) }}, 
\end{align*}
where the last inequality holds from Jensen's inequality, because $1/\sqrt{1 + a x}$ is a convex function of $x$ for any constant $a > 0$ and $\lt(\frac{n p_w(1 - p_w)^{n}}{\ee N(w)}, n \geq 0\rt)$ is a probability distribution.
Furthermore, using $\ee N(w)^2 = \ee N(w)\frac{(2 - p_w)}{p_w} \leq 2\,\ee N(w)/p_w$, we write
$\ee \wt N(w) \geq \ee N(w)\frac{1}{\sqrt{1 + 2 d_w z_1^2/p_w}}.$
Since the distribution \eqref{eqn:time_var_min_pmf} is not the variance minimizing distribution, we have
\begin{align*}
\ee\wt N(w)\var\, \wt \beta(w) &\geq \ee \wt N(w) \var\,\whbeta(w) \geq \ee N(w) \var\, \whbeta(w)\lt( \frac{1}{\sqrt{1 + 2 d_w z_1^2/p_w}}\rt).
\end{align*}
We establish \eqref{eqn:tvp_order} by showing that $\frac{1}{\sqrt{1 + 2 d_w z_1^2/p_w}} = 1 + O(w)$. From the definition of $d_w$,
\begin{align*}
\frac{w z_1}{p_w}\sum_{n = 0}^\infty \frac{p_w(1- p_w)^n}{\sqrt{1 + d_w z_1^2\, n}} = 1,
\end{align*}
and thus, using $p_w \leq w z_1$, we write that
$1 - w\, z_1 \leq (1 - p_w)\sum_{n = 1}^\infty \frac{p_w(1- p_w)^{n-1}}{\sqrt{1 + d_w\, z_1^2\, n}} \leq (1 - p_w)\frac{1}{\sqrt{1 + d_w\, z_1^2}}.$
Consequently, $1 + d_w\, z_1^2 \leq \lt(\frac{1 - p_w}{1 - w\, z_1}\rt)^2$. Hence, using \eqref{eqn:pw_bnd}, $1 + d_w\, z_1^2   = 1 + O(w^2)$, that is,
$d_w = O(w^2)$ and thus, $d_w/p_w = O(w)$ because $p_w = O(w)$. This concludes that $\frac{1}{\sqrt{1 + 2 d_w z_1^2/p_w}} = 1 + O(w)$ and hence establishes~\eqref{eqn:tvp_order}.

\subsection{Proof of Theorem~\ref{thm:impl}}
\label{sec:proof_impl}
From the assumption, we have $w_k \leq \frac{M_1(k)}{M_2(k)}$. Therefore, similar to \eqref{eqn:pw_bnd}, we can obtain that 
\begin{align}
\label{eqn:bnd_P_k}
0 < w_k\, M_1(k) - w_k^2\, M_2(k)/2 \leq P_k \leq w_k\, M_1(k).
\end{align}
From the definitions of $w_k, M_1(k)$ and $M_2(k)$, it is easy to see that $\lim_{k \to \infty} kw_k\, M_1(k) = 1,\, a.s.$.
Using the upper bound in \eqref{eqn:bnd_P_k}, we show that $\delta_k
:= 1 - \frac{ 1 + w^2_k\, z_2 - 2 w_k z_1}{1 - P_k}$ is lower bounded
by $\frac{w_k}{1- P_k} \lt(2z_1 - M_1(k) - w_k\, z_2\rt)$.  Since
$\lim_{k \to \infty} M_1(k) = z_1, \, a.s.$ and $\lim_{k \to \infty}
kw_k = 1/z_1,\, a.s.$, for every realization of $\wt Z_1, \wt Z_2,
\dots$, there exists a $K$ such that $\delta_k > 0$ for all $k \geq
K$, and hence $\var\,\lt( \whbeta(w_k) \big|M_1(k), M_2(k)\rt)$ is finite
and equal to $\frac{w_k^2}{P_k\delta_k} - \beta^2$. It is now enough
to show that
\begin{align*}
k\, \lt( \frac{w_k^2\, z_1^2}{P_k\delta_k} - 1 \rt) \longrightarrow \frac{\var\,Z}{z_1^2}  \lt[1 + \chi^2_1 \rt],\, a.s., \quad \text{as } k \to \infty.
\end{align*}
Write $\displaystyle \lt( \frac{w_k^2\, z_1^2}{P_k\delta_k} - 1 \rt) =
\lt(\frac{w_k^2\, z_1^2}{p^2_{w_k}} - 1 \rt) + w_k^2\,
z_1^2\lt(\frac{1}{P_k\delta_k} - \frac{1}{p^2_{w_k}} \rt)$, where
$p_{w_k} = 1 - \sqrt{1 + w^2_k\, z_2 - 2w_k\, z_1}$. Using
\eqref{eqn:cost_var} and $\lim_{k \to \infty} k\, w_k = 1/z_1$, we
have $\lim_{k \to \infty}k\lt(\frac{w_k^2\, z_1^2}{p^2_{w_k}} - 1 \rt)
= \frac{\var\, Z}{z_1^2},\, a.s$.  By simplifying the terms in
$P_k\delta_k - p^2_{w_k}$, we have $\frac{1}{P_k\delta_k} -
\frac{1}{p^2_{w_k}} = \frac{1}{P_k\delta_k (1- P_k)} \lt(\frac{P_k -
  p_{w_k}}{p_{w_k}}\rt)^2.$ Since $P_k -p_{w_k} = w_k\,(M_1(k) - z_1)
(1+ O(w_k))$, we can write $\frac{P_k - p_{w_k}}{p_{w_k}} =
\frac{w_k}{p_{w_k}} (M_1(k) - z_1)(1 + O(w))$. Therefore, using
$\lim_{k \to \infty} \frac{w_k^2\, z_1^2}{P_k\delta_k(1 - P_k)} = 1,\,
a.s.$ and the fact that asymptotically $\sqrt{k}(M_1(k) - z_1)$ has
a zero-mean normal distribution with variance $\var\, Z$, we complete the
proof with the observation that $w_k/p_{w_k} \to 1/z_1$ as $k \to
\infty$.

\section*{Acknowledgements}
The work of the first and the second authors has been supported by the Australian Research Council Centre of Excellence for Mathematical and Statistical Frontiers (ACEMS), under grant number CE140100049.  
\bibliographystyle{apalike}
\bibliography{Ref}
\end{document}